\documentclass[12pt]{article}

\usepackage[T2A]{fontenc}
\usepackage[utf8]{inputenc}

\usepackage{amsfonts}
\usepackage{amssymb}
\usepackage{amsmath}
\usepackage{amsthm}

\def \le {\leqslant}
\def \ge {\geqslant}

\begin{document}

\begin{Large}
 \centerline{\bf Sur une question de N. Chevallier liée à   }
 \centerline{\bf l'approximation Diophantienne simultanée}
\end{Large}
\vskip+0.5cm

\centerline{ par {\bf Nikolay Moshchevitin}\footnote{la recherche est financée par la subvention de RFBR No.12-01-00681-a
et par la subvention de le Gouvernement Russe, projet 11. G34.31.0053.
}}

\vskip+0.5cm
\begin{small}
 {\bf Résumé.}\, Nous prouvons une conjecture proposée par Nicolas Chevallier qui concerne des matrices unimodulaires liées à l'approximation Diophantienne simultanée des nombres réels.
\end{small}
\vskip+0.5cm

{\bf 1. Approximation Diophantienne simultanée.}
  
Soit $n$ un entier naturel. Dans cet article nous considérons un vecteur réel $\xi$ de la forme 
$\xi = (1,\xi_1,...,\xi_n)$ dont les coordonnées sont linéairement indépendants sur $\mathbb{Z}$.
Nous sommes intéressés à l'approximation du sous-espace vectoriel engendré par $\xi$ par des points entiers ${\bf z} =(q,a_1,...,a_n)$.
On considère la fonction
$$
\psi_\xi (t) = \min_{q\in \mathbb{Z}, 1\le q \le t} \, \max_{1\le k \le n} ||q\xi_k||,
$$
où $||\cdot ||$ désigne la distance à l'entier le plus proche. Cette fonction est décroissante et constante par morceaux.
Soient
$$
q_0 =1\le q_1< q_2 <....<q_\nu<q_{\nu+1}
<...
$$
les sauts de $\psi_\xi (t)$. Ils alors correspondent aux {\it vecteurs meilleures approximations} ${\bf g}_\nu = (q_\nu, a_{1,\nu},...,a_{n,\nu}) \in \mathbb{Z}^{n+1}$ qui sont définis par les conditions
$$
|| q_\nu \xi_k||
= |q_\nu \xi_k - a_{k,\nu}|,\,\,\,
1\le k \le n.
$$
Notons que d'après le théorème de Minkowski sur les corps convexes on sait que
$$
 \psi_\xi (t)\le t^{-\frac{1}{n}},
$$
ou
\begin{equation}\label{apo}
\max_{1\le k \le n} ||q_\nu\xi_k|| \le q_{\nu+1}^{-\frac{1}{n}}.
\end{equation}
Dans le cas $ n = 1$ d'après la théorie de fractions continues nous savons que pour les approximations (\ref{apo}) nous avons une borne inférieure du même ordre, c'est
$$
(2q_{\nu+1})^{-1}<
||q_\nu \xi_1 || <q_{\nu+1}^{-1}.
$$
Aussi nous savons que
$$
\left|
\begin{array}{cc}
q_\nu & q_{\nu+1}\cr
a_{\nu,1}&a_{\nu+1,1}
\end{array}
\right| = \pm 1.
$$
Ces simples observations conduisent au corollaire suivant.

{\bf Proposition 1.}\,\,{\it
Pour nombre irrationnel $\xi_1$ il existe une infinité de matrices unimodulaires
$$
\left(
\begin{array}{cc}
q' & q''\cr
a_{1}'&a_{1}''
\end{array}
\right)
$$
telles que
$$
\max\left\{
q'{|q'\xi_1-a_1'|},q''{|q''\xi_1-a_1''|}
\right\}
\le 1.
$$
}
La situation dans le cas $n \ge 2$ est tout à fait différent.
Considérons une fonction décroissante $\varphi (t)$ telle que
\begin{equation}\label{tech}
\varphi (t) =o(1), \,\,\, t\to\infty.
\end{equation}
Le résultat principal de cet article est le théorème suivant.

{\bf Théorème 1.} \,{\it Pour une foncion donnée $\varphi (t)$ qui deminue à zéro lorsque $ t\to \infty$, il existe un vecteur $(1,\xi_1,\xi_2) \in \mathbb{R}^3$ dont les composants sont linéairement indépendants sur $\mathbb{Z}$ et tel que pour toute matrice entière 
\begin{equation}\label{meer}
M=
\left(
\begin{array}{ccc}
q'&q''&q'''\cr
a_{1}'& a_{1}''&a_{1}'''\cr
a_{2}'& a_{2}''&a_{2}'''
\end{array}
\right),\,\,\,\, {\rm det M} =\pm 1,\,\,\,\,\,  q',q'',q'''\ge 1
\end{equation}
on a
$$
\max \left\{ \frac{\max_{j=1,2} |q'\xi_j - a_j'|}{\varphi (q')},
 \frac{\max_{j=1,2} |q''\xi_j - a_j''|}{\varphi (q'')},
 \frac{\max_{j=1,2} |q'''\xi_j - a_j'''|}{\varphi (q''')}\right\}
\ge \varepsilon,
$$
avec un certain nombre positif $\varepsilon $ qui dépend de la fonction $\varphi$.
}

Notre Théorème 1 donne une réponse affirmative à la question posée par N. Chevallier dans \cite{che} dans le cas $s=2$.
Bien sûr, un résultat similaire sera vrai pour tout $ n \ge 2$.

D'après l'argument ci-dessous, il est possible de prouver la déclaration suivante.

{\bf Proposition 2.}
{\it 
Il existe $\xi_1, \xi_2$ linéairement indépendants sur $\mathbb{Z}$ avec $1$ et tels que le déterminant
d'une matrice de la forme
$$
\left(
\begin{array}{ccc}
q_{\nu_1}&q_{\nu_2}&q_{\nu_3}\cr
a_{1,\nu_1}& a_{1,\nu_2}&a_{1,\nu_3}\cr
a_{2,\nu_1}& a_{2,\nu_2}&a_{2,\nu_3}
\end{array}
\right),\,\,\,\,\, \nu_1<\nu_2<\nu_3
$$
(ici $(q_\nu, a_{\nu,1},a_{\nu,2})$ sont les vecteurs meilleures approximations)
n'est jamais égal à $\pm 1$. }

Bien sûr, ce résultat est aussi valable dans toute dimension $ n \ge 2$.

Nous tenons à rappeler deux résultats liés à l'approximation Diophantienne simultanée.

Le premier résultat va jusqu'à V. Jarn\'{\i}k \cite{ja} (voir aussi \cite{D} et  \cite{L}). Il déclare que pour $ n \ge 2$ il y a un nombre infini de triplets
de {\it linéairement indépendants} vecteurs meilleures approximations consécutifs 
${\bf g}_{\nu}, {\bf g}_{\nu+1}, {\bf g}_{\nu+2}$.

La seconde est due à l'auteur \cite{UMN}. Il déclare que pour tout $n\ge 2$ il existe $\xi_1,...,\xi_n$ linéairement indépendants sur $\mathbb{Z}$ avec $1$  et tels que pour tout $\nu$
la matrice
$$
\left(
\begin{array}{cccc}
q_{\nu}&q_{\nu +1}&...&q_{\nu+n}\cr
a_{1,\nu}& a_{1,\nu+1}&...&a_{1,\nu+n}\cr
....& ....&...&...\cr
a_{n,\nu}& a_{n,\nu+1}&...&a_{n,\nu+n}
\end{array}
\right)
$$
de $n+1$ vecteurs meilleures approximations consécutifs est de rang $\le 3$.
Ce résultat donne un contre-exemple à la conjecture de Lagarias \cite{L}.

Pour plus d'informations concernant les vecteurs meilleures approximations nous nous référons à \cite{che, M1,M2}

Pour conclure cette section, nous formulons un résultat du type de Jarn\'{\i}k qui est basé sur une observation simple par V.  Jarn\'{\i}k
(voir \cite{TB}, Satz  9 et Théorème 17 de \cite{M2}). 
Pour une matrice $M$ de la forme (\ref{meer}) nous définissons
$$
R(M) =
\max \left\{ {\max_{j=1,2} |q'\xi_j - a_j'|},
{\max_{j=1,2} |q''\xi_j - a_j''|},
{\max_{j=1,2} |q'''\xi_j - a_j'''|}
\right\}.
$$

{\bf Théorème 2.} \,{\it 
Supposons que $1,\xi_1,\xi_2$ sont linéairement indépendants sur $\mathbb{Z}$. Il existe la suite $M_\nu$ des matrices de la forme (\ref{meer}) telle que
$R(M_\nu) \to 0$ lorsque $\nu \to \infty$.}
\vskip+0.2cm

Théorème 1 sera prouvée dans les sections 2 - 5.
Nous donnons une preuve du Théorème 2 dans la section 6.

{\bf 2. Lemmes.}

Ici $|\cdot |$ signifie la norme Euclidienne,
$ {\rm dist} \, ({\cal A}, {\cal B})$ désigne la distance Euclidienne entre les ensembles ${\cal A}, {\cal B}$.
Par ${\rm angle} \, ({\bf u}, {\bf v})$
nous notons l'angle entre
les vecteurs ${\bf u}$ et ${\bf v}$.
Par ${\rm angle} \, (L, P)$
nous notons aussi l'angle entre les sous-espaces $L$ et   $P$ de dimension un ou deux.

Pour deux vecteurs indépendants ${\bf z}', {\bf z}''\in \mathbb{R}^3$ nous définissons le sous-espace vectoriel
$$
L({\bf z}',{\bf z}'') = {\rm span} ({\bf z}',{\bf z}').
$$
Pour ${\bf z}', {\bf z}''\in \mathbb{Z}^3$ nous considérons le réseau
$$
\Lambda ({\bf z}',{\bf z}'') = L({\bf z}',{\bf z}'') \cap \mathbb{Z}^3.
$$
Notons que dans le cas où une paire 
${\bf z}', {\bf z}''\in \mathbb{Z}^3$
peut être complétée à une base de tout $\mathbb{Z}^3$ on a
$$
\Lambda ({\bf z}',{\bf z}'') = \langle{\bf z}',{\bf z}''\rangle_\mathbb{Z}.
$$ 
Pour deux points entiers ${\bf z}', {\bf z}''$ sous la condition
\begin{equation}\label{root}
L( {\bf z}', {\bf z}'') \cap \mathbb{Z}^3 = \langle{\bf z}', \bf {z}''\rangle_\mathbb{Z}
\end{equation}
nous considérons un point $ {\bf y} ({\bf z}', {\bf z}'')$ qui complète la paire ${\bf z}', {\bf z}''$
à une base de $\mathbb{Z}^3$.
Ensuite, nous définissons deux sous-espaces affines de dimension deux
\begin{equation}\label{recc}
L^{\pm} ({\bf z}', {\bf z}'' ) =
L({\bf z}',{\bf z}'')\pm  {\bf y} ({\bf z}', {\bf z}'').
\end{equation}
Nous devons d'introduire une quantité plus. Pour deux points entiers ${\bf z}', {\bf z}''$ sous la condition (\ref{root}) nous considérons la valeur
\begin{equation}\label{vaal}
\eta ({\bf z}', {\bf z}'') = \min_{{\bf x}\in \Lambda  ({\bf z}',{\bf z}'')   \setminus {\rm span}({\bf z}')}   {\rm dist} ({\bf x}, {\rm span} ({\bf z}'))> 0
.
\end{equation}

{\bf Lemme 1.}\,\,{\it
Soit ${\bf z} \in \mathbb{Z}^3$. Supposons que ${\bf z}', {\bf z}''\in \mathbb{Z}^3$  
satisfait (\ref{root}) et $ {\bf z} \not\in  L( {\bf z}', {\bf z}'') $.
Considérons le point ${\bf w} = L^+ ({\bf z}', {\bf z}'' )\cap {\rm span} \, ({\bf z})$. Soit $\delta$ l'angle
$$
  {\rm angle}\, (  L( {\bf z}', {\bf z}''), {\rm span} \, ({\bf z}))>\delta >0.
$$
Supposons que
$$
{\bf x}\in  L^- ({\bf z}', {\bf z}'' )
\cup L^+ ({\bf z}', {\bf z}'' )
$$
et
\begin{equation}\label{aar}
|{\bf x}|\ge 2|{\bf w}|  .
\end{equation}
Alors
$$
{\rm dist }\, ({\bf x},  {\rm span} \, ({\bf z})) \ge  \frac{|{\bf x}|}{2}   \sin \delta .
$$
}

Preuve.\,\, 
Notons par ${\bf x}^*$ la projection orthogonale du point ${\bf x}$
sur le sous-espace  ${\rm span} \, ({\bf z})$, qui est de dimension un.  Alors,
$$
{\rm dist }\, ({\bf x},  {\rm span} \, ({\bf z})) =
{\rm dist }\, ({\bf x},   {\bf x}^*) =
|{\bf x} \pm {\bf w}| \sin \theta,
$$
où $\theta \ge \delta$ est l'angle entre  ${\bf x} \pm {\bf w}$ et ${\bf z}$.
Ainsi
$$
{\rm dist }\, ({\bf x},  {\rm span} \, ({\bf z}))
\ge (|{\bf x}|-|{\bf w}|)  \sin  \delta \ge \frac{|{\bf x}|}{2}   \sin \delta 
 $$
(ici nous utilisons (\ref{aar})).
Lemme 1 est prouvé.$\Box$

\vskip+0.2cm

D'après Lemme 1 nous déduisons immédiatement

{\bf Corollaire 1.}\,\,{\it
Soit ${\bf z} \in \mathbb{Z}^3$. Supposons que ${\bf z}', {\bf z}''\in \mathbb{Z}^3$  
satisfait (\ref{root}) et $ {\bf z} \not\in  L( {\bf z}', {\bf z}'') $.
Alors il existe des nombres positifs $\delta ({\bf z}, {\bf z}', {\bf z}'')$ et $T ({\bf z}, {\bf z}', {\bf z}'')$
tels que pour tous les vecteurs $\xi = (1,\xi_1,\xi_2)$ sous la condition
$$
{\rm angle}\, (\xi , {\bf z}) \le   \delta ({\bf z}, {\bf z}', {\bf z}'')
$$
et pour tous les vecteurs entiers ${\bf x} = (q, a_1, a_2), q\ge 1$ sous la condition
 $$
{\bf x}\in  L^- ({\bf z}', {\bf z}'' )
\cup L^+ ({\bf z}', {\bf z}'' )
\,\,\,\text{
et}
\,\,\,
|{\bf x}|\ge T ({\bf z}, {\bf z}', {\bf z}'')
$$
on a
\begin{equation}\label{result}
\max_{j=1,2} |q\xi_j - a_j|\ge \varphi (q).
\end{equation}
}
Il est clair que dans Corollaire 1 on peut considérer une collection finie de couples ${\bf z}', {\bf z}''$.
Donc nous avons la déclaration suivante.

{\bf Corollaire 2.}\,\,{\it
Soit ${\bf z} \in \mathbb{Z}^3$.
Soit $\hbox{\got C}$ une collection finie de couples $({\bf z}', {\bf z}'')$ de points entiers tels que chacun d'eux
satisfait (\ref{root}) et $ {\bf z} \not\in  L( {\bf z}', {\bf z}'') $.
Alors il existe des nombres positifs $\delta ({\bf z}, \hbox{\got C})$ et $T ({\bf z}, \hbox{\got C})$
tels que pour tous les vecteurs $\xi = (1,\xi_1,\xi_2)$ sous la condition
$$
{\rm angle}\, (\xi , {\bf z}) \le   \delta ({\bf z},\hbox{\got C})
$$
et pour tous les vecteurs entiers ${\bf x} = (q, a_1, a_2), q\ge 1$ sous la condition
 $$
{\bf x}\in
\bigcup_{ 
({\bf z}', {\bf z}'')\in \hbox{\got C}}
( L^- ({\bf z}', {\bf z}'' )
\cup L^+ ({\bf z}', {\bf z}'' )
)
\,\,\,\text{
et}
\,\,\,
|{\bf x}|\ge T ({\bf z}, \hbox{\got C})
$$
on a (\ref{result}).
 }

\vskip+0.2cm

Considérons un point entier $ {\bf z} \in \mathbb{Z}^3 \setminus \{{\bf 0}\}$.
Soit $\Lambda\subset \mathbb{Z}^3$ un sous-réseau de dimension deux tel que $\Lambda\ni {\bf z}$ et
$\varepsilon >0$, nous considérons l'ensemble des sous-réseaux
\begin{equation}\label{sub}
{\cal L} ={\cal L} ({\bf z},\Lambda,\varepsilon) =\{
\Lambda '\subset \mathbb{Z}^3: \,\, {\rm dim }\,\Lambda ' = 2,\,\,   {\bf z } \in \Lambda',\,\,{\rm angle} \, ( {\rm span}\,\Lambda ', {\rm span}\, \Lambda)<\varepsilon\}.
\end{equation}

{\bf Lemme 2.}\,\,{\it
Considérons un ensemble
${\cal L} ({\bf z},\Lambda,\varepsilon) $
de la forme (\ref{sub}). Soit $T$ un nombre positif. Alors, il existe un réseau $\Lambda \in {\cal L} ({\bf z},\Lambda,\varepsilon) $
tel que
pour tout point ${\bf x}$ qui satisfait
$$
{\bf x} \in \Lambda\,\,\,\text{et}\,\,\, |{\bf x} |\le T
$$
on a
$$
{\bf x }\in {\rm span}\,( {\bf z}).
$$
}
Preuve. \,\,Le lemme résulte de l'observation que tout ensemble de la forme (\ref{sub}) contient un nombre infini d'éléments.$\Box$

\vskip+0.2cm
{\bf Lemme 3.}\,\,{\it Soit ${\bf z}$ un point entier et
$\Lambda\ni {\bf z}$ un sous-réseau de dimension deux de $\mathbb{Z}^3$.
Considérons un point entier ${\bf z}'$ indépendant de ${\bf z}$.
Alors il existe positif $\varepsilon^*$ et un sous-réseau de dimension deux $ \Lambda_*\ni {\bf z}$ tels que
 $${\cal L}_*= {\cal L} ({\bf z}, \Lambda_*, \varepsilon_*)\subset {\cal L} ({\bf z},\Lambda,\varepsilon) $$
et
pour tout 
$
\Lambda \in {\cal L}_*$ on a
$$\Lambda \cap (L^-({\bf z}, {\bf z}') \cup L^+ ({\bf z}, {\bf z}') ) = \varnothing
.
 $$
 }

Preuve.\,\,
Le sous-réseau affine $L^-({\bf z}, {\bf z}')\cap\mathbb{Z}^3$ (et le sous-réseau $L^+({\bf z}, {\bf z}')\cap\mathbb{Z}^3$) se divise en sous-réseaux affines (de dimension un) $\Gamma_i$ qui est parallèle à
${\rm span}, ({\bf z})$:
$$
L^-({\bf z}, {\bf z}')\cap\mathbb{Z}^3 = \bigsqcup_{i \in \mathbb{Z}} \Gamma_i.
$$
Il suffit de traiter avec $L^-({\bf z}, {\bf z}')$, par l'argument de la symétrie.
Nous considérons deux points différents $ {\bf w}_j\in L^-({\bf z}, {\bf z}'), j=1,2$
tels qu'ils
appartiennent à deux sous-espaces affines voisins (de dimension un) $\Gamma_i$ et $\Gamma_{i+1}$, respectivement,
et pour certains ${\bf w}$ de l'intervalle ouvert avec les bornes ${\bf w}_1, {\bf w}_2$ le sous-espace $ L ({\bf w}, {\bf z})$ contient un sous-réseau
$$
\Lambda_*  = L ({\bf w}, {\bf z})  \cap \mathbb{Z}^3 \in {\cal L}.
$$
Il est clair que pour un certain petit
$\varepsilon_*$ l'ensemble ${\cal L}_*$ de la forme (\ref{sub}) satisfait la propriété désirée.$\Box$

\vskip+0.2cm

D'après Lemme 3 on déduit immédiatement le suivant

{\bf Corollaire 3.}\,\,{\it
Soit $\hbox{\got E}$ une collection finie de couples ${\bf z}'$ et chacun de ces points est indépendant de ${\bf z}$.
Supposons que $ \Lambda \ni {\bf z}$ est un sous-réseau entier de dimension deux. Alors il existe un ensemble
$${\cal L}_* ={\cal L} ({\bf z}, \Lambda_*, \varepsilon_*)$$
de la forme (\ref{sub}) tel que pour tout $\Lambda \in {\cal L}_*$ on a
$$\Lambda \cap \left( \bigcup_{{\bf z}'\in \hbox{\got E}} (L^-({\bf z}, {\bf z}') \cup L^+ ({\bf z}, {\bf z}') ) \right)= \varnothing
.
 $$

}

\vskip+0.2cm

Pour deux points indépendants ${\bf z}$ et ${\bf z}^*$  nous considérons l'ensemble
$$
{\cal P}({\bf z},{\bf z}^*)
=
\bigcup_P P,
$$
où l'union est prise sur tous les sous-espaces vectoriels de dimension deux $P$
tels que
$$
{\rm span}\, ({\bf z}) \subset P\,\,\,\text{et}\,\,\, {\rm angle} (P, L({\bf z},{\bf z}^*))\ge \frac{3\pi}{8}.
$$
Nous avons aussi besoin d'un ensemble
$$
\overline{\cal P}({\bf z},{\bf z}^*)
=
\bigcup_P P,
$$
où l'union est prise sur tous les sous-espaces linéaires de dimension deux $P$
tels que
$$
{\rm span}\, ({\bf z}) \subset P\,\,\,\text{et}\,\,\, {\rm angle} (P, L({\bf z},{\bf z}^*))\ge \frac{\pi}{4}.
$$
Il est clair que
$$
\overline{\cal P}({\bf z},{\bf z}^*)
\supset
{\cal P}({\bf z},{\bf z}^*).
$$

{\bf Lemme 4.}\,\, {\it
Supposons que les points entiers ${\bf z}$ et ${\bf z}^*$ peuvent être complétés à une base de $\mathbb{Z}^3$
Alors
$$
{\rm dist} (\overline{\cal P}({\bf z},{\bf z}^*), \Lambda ({\bf z},{\bf z}^*)\setminus {\rm span}({\bf z}) )  \ge
\frac {\eta({\bf z},{\bf z}^*)}{\sqrt{2}},
$$
où $\eta(\cdot,\cdot )$ est défini dans (\ref{vaal}).
}
 
Preuve.\,\, La distance
entre $ {\bf x}\in \Lambda ({\bf z},{\bf z}^*)\setminus {\rm span}({\bf z})$ et $\overline{\cal P}({\bf z},{\bf z}^*)$ n'est pas inférieure à la distance entre
 $ {\bf x}$ et ${\rm span}({\bf z})$ multiplié par $\sin \frac{\pi}{4}$.$\Box$

\vskip+0.2cm
{\bf 3. Vecteurs ${\bf z}_\nu$.}

Dans cette section nous construisons une suite de vecteurs d'entiers ${\bf z}_\nu$ par une certaine procédure inductive.
On met
$$
{\bf z}_1 = (1,0,0),\,\,\,\,\,\,{\bf z}_2 = (0,1,0).
$$
Maintenant, nous supposons que les vecteurs $ {\bf z}_1,....,{\bf z}_\nu, \nu\ge 2$ sont déjà définis.

Pour une fonction décroissante $\varphi  (t)$  nous définissons la fonction
 $\phi (t)$ qui est la fonction inverse de $\varphi (t)$.
Nous définissons
\begin{equation}\label{hanu1}
H_\nu = \phi (2^{-3}\eta ( {\bf z}_{\nu-1}, {\bf z}_\nu)).
\end{equation}
Considérons les ensembles
 $$
\hbox{\got C}_\nu  =
\bigcup_{\lambda=1}^\nu
 \bigcup_{\mu=1}^{\nu-1}
\{ ({\bf z}', {\bf z}'') \,\,\text{satisfait (\ref{root})}:\, {\bf z}_\nu\not\in L({\bf z}', {\bf z}''),\,\,
{\bf z}' \in \Lambda_\lambda,\, |{\bf z}'| \le H_\lambda,\, 
{\bf z}'' \in \Lambda_\mu  ,\, |{\bf z}''| \le H_\mu\}
$$
et
$$
\hbox{\got E}_\nu 
=\bigcup_{\mu=1}^{\nu}
\{ {\bf z}':  ({\bf z}_\nu, {\bf z}') \,\,\text{satisfait (\ref{root})},
\,{\bf z}' \in \Lambda_\mu \, |{\bf z}'| \le H_\mu\}
$$
Il est claire que $
\hbox{\got C}_\nu 
$
et
$
\hbox{\got E}_\nu 
$
sont des ensembles finis.

Mettons
$$
\delta_\nu = \delta ({\bf z}_\nu, \hbox{\got C}_\nu),\,\,\,\,
T_\nu = T ({\bf z}_\nu, \hbox{\got C}_\nu)
,\,\,\,\, \nu \ge 2,
$$
où $\delta(\cdot,\cdot)$ et $T(\cdot,\cdot )$ sont définis dans Corollaire 2. Bien sûr, nous pouvons supposer que $ \delta_\nu <\delta_{\nu-1}/2$,
où $\delta_{\nu-1}$ est défini à l'étape précédente de la construction
(au début du processus on met
$\delta_1=\delta_2 = \pi, T_1 = T_2 = 1$).

Nous
supposons que
\begin{equation}\label{daba}
{\rm angle} ({\bf z}_\nu, {\bf z}_{\nu-1}) <\frac{\delta_{\nu-1}}{2}.
\end{equation}

Nous devons définir maintenant le vecteur ${\bf z}_{\nu+1} $ de telle manière que le couple $ ({\bf z}_\nu, {\bf z}_{\nu+1})$ peut être complétée à une base de
$\mathbb{Z}^3$ et définir le réseau correspondant
$
\Lambda_{\nu+1} = \langle {\bf z}_\nu, {\bf z}_{\nu+1}\rangle_\mathbb{Z}.
$
Nous allons le faire de la manière suivante.
{\it Au début} nous allons définir un réseau $
\Lambda_{\nu+1}\ni {\bf z}_\nu$ et {\it ensuite} nous allons choisir le vecteur $ {\bf z}_{\nu+1}$ pour compléter le couple $ ({\bf z}_\nu, {\bf z}_{\nu+1})$.

Nous prenons un réseau  $\Lambda_{\nu+1}$ pour satisfaire les conditions

({\bf i})  $ {\bf z}_\nu \in \Lambda_{\nu+1}$,

({\bf ii}) $\mathbb{Z}^3 \cap {\rm span } \Lambda_{\nu+1} = \Lambda_{\nu+1}$,

({\bf iii})  $\Lambda_{\nu+1} \subset {\cal P}({\bf z}_\nu, {\bf z}_{\nu-1})$,

({\bf iv}) pour chaque $ {\bf x} \in \Lambda_{\nu+1}$ tel que $ |{\bf x}|\le T_\nu$ on a $ {\bf x}\in {\rm span} ({\bf z}_\nu)$,

({\bf v})   $\Lambda_{\nu+1} \cap \left( \bigcup_{{\bf z}'\in \hbox{\got E}_\nu} (L^-({\bf z}, {\bf z}') \cup L^+ ({\bf z}, {\bf z}') ) \right)= \varnothing
.
 $

L'existence d'un tel réseau $\Lambda_{\nu+1}$ résulte de Lemme 2 and Corollaire 3.

Maintenant, nous expliquons comment choisir le point entier $ {\bf z}_{\nu+1}=(q_{\nu+1},a_{1,\nu+1}, a_{2,\nu+1})\in \Lambda_{\nu+1}$.
Comme $ {\bf z}_{\nu} \in \Lambda_{\nu+1}$ est un point primitif,
le réseau $\Lambda_{\nu+1}$ se divise en réseaux affines parallèles (de dimension un) $\Gamma^k$ de telle manière que
$
\Lambda_{\nu+1}=\bigcup_{k\in \mathbb{Z}}\Gamma^k$, $ \Gamma^0 = {\rm span } ({\bf z}_{\nu}) \, \cap \mathbb{Z}^3$. Ici $\Gamma^{\pm1}$ sont les réseaux les plus proches de $\Gamma^0$.
Si nous prenons ${\bf z}_{\nu+1} \in \Gamma^1$, nous voyons que le couple
$({\bf z}_\nu, {\bf z}_{\nu+1})$ peut être complétée à une base de $\mathbb{Z}^3$. Notons que si $|{\bf z}_{\nu+1}|$ est assez grand
(et donc $q_{\nu+1}$ est grand) alors l'angle $ {\rm angle} ( {\bf z}_{\nu+1}, {\bf z}_\nu)$ est petit.
Il est clair que cet angle tend vers zéro lorsque $|{\bf z}_{\nu+1}|$  tend vers l'infini.
Il existe donc
$$
W_\nu^1 = W_\nu (\delta ({\bf z}_\nu, \hbox{\got C}_\nu))
$$
tel que
si  $|{\bf z}_{\nu+1}|\ge W_\nu^1 $ alors
\begin{equation}\label{delta1}
 {\rm angle} ( {\bf z}_{\nu+1}, {\bf z}_\nu)<\delta_\nu/2.
\end{equation}
On peut supposer que $ \delta_\nu$ est assez petit, alors
\begin{equation}\label{delta2}
\{ {\bf x}\in \mathbb{R}^3: \,
 {\rm angle} ( {\bf x}, {\bf z}_\nu)<\delta_\nu\} \subset
\{ {\bf x}\in \mathbb{R}^3: \,
 {\rm angle} ( {\bf x}, {\bf z}_{\nu-1})<\delta_{\nu-1}\}
.
\end{equation}

Si nous choisissons ${\bf z}_{\nu+1}$ nous pouvons considérer la distance Euclidienne
$$
\rho_\nu = 
{\rm dist}\, ({\bf z}_\nu, {\rm span}({\bf z}_{\nu+1})).
$$
Nous voyons d'après la définition que $\rho_\nu$ dépend du choix du point ${\bf z}_{\nu+1}$.
Notons que pour tout choix de ${\bf z}_{\nu+1} \in \Gamma^{\pm1}$ la valeur
$q_{\nu+1}\rho_\nu$ sera du même ordre que le volume fondamental du réseau de dimension deux $\Lambda_{\nu+1}$, qui est déjà définie.
Plus précisément, la quantité
$\frac{q_{\nu+1}\rho_\nu}{{\rm det}\Lambda_{\nu+1}}
$
est bornée par l'arrière et est séparé de zéro par une constante positive.

Mais $ \rho_{\nu}$ tend vers zéro lorsque $|{\bf z}_{\nu+1}|$ tend vers l'infini.
Il existe donc
$$
W_\nu^2 = W_\nu ( {\bf z}_{\nu-1}, {\bf z}_\nu)
$$
tel que pour tous ${\bf z}_{\nu+1}\in \Gamma^{\pm1}$ sous la condition  $|{\bf z}_{\nu+1}|>W_\nu^2$
on a
$\rho_{\nu} \le \rho_{\nu-1}/2$ et
\begin{equation}\label{ce}
 \varphi \left(\frac{1}{64\rho_{\nu-1}\rho_\nu}\right) \le \frac{1}{16q_{\nu+1} \rho_{\nu}}
\end{equation}
(ici $\rho_{\nu-1}$ supposé être défini par des moyens des points de ${\bf z}_{\nu-1}, {\bf z}_\nu$ à l'étape précédente de la construction).

Maintenant, nous fixons ${\bf z}_{\nu+1} $ avec
$$
|{\bf z}_{\nu+1}|>\max ( W_\nu^1, W_\nu^2)
.
$$
Nous avons construit le point suivant ${\bf z}_{\nu+1}$. Pour le point construit les conditions (\ref{delta1},\ref{delta2},\ref{ce}) sont valables.
Et le réseau correspondant $\Lambda_{\nu+1}$ satisfait ({\bf i}) - ({\bf v}).
De plus, notre construction donne (\ref{daba}) avec $\nu$ remplacé par $\nu+1$.

Maintenant, nous mettons
$$
\xi_{j, \nu} = \frac{a_{i,\nu}}{q_\nu},\,\,\,j = 1,2
$$
et
\begin{equation}\label{xi}
\xi_j = \lim_{\nu\to \infty} \xi_{j,\nu},\,\,\,\
\xi = (1,\xi_1,\xi_2).
\end{equation}
Bien sûr, nous pouvons supposer que $ \xi_1,\xi_2 \in (0,1)$.

Nous voyons d'après (\ref{delta1},\ref{delta2}) que
\begin{equation}\label{dee}
\xi \in 
\{ {\bf x}\in \mathbb{R}^3: \,
 {\rm angle} ( {\bf x}, {\bf z}_\nu)<\delta_\nu\} \,\,\,\,\,
\forall \nu.
\end{equation}
Notons que pour $\xi$ défini dans (\ref{xi})
nous avons
$$
{\rm dist }( {\bf z}_\nu, {\rm span }(\xi))
\le \sum_{k =\nu}^\infty {\rm dist} ({\bf z}_k, {\rm span} ({\bf z}_{k+1})) \le 2\rho_\nu,
$$
et ainsi
\begin{equation}\label{uro}
\max_{j=1,2} |q_\nu \xi_j - a_{\nu, j} | \le 4\rho_\nu.
\end{equation}
Par ailleurs il faut noter que d'après {\bf (iii)}, il en résulte que
$$
\xi \in  \overline{\cal P}({\bf z}_\nu, {\bf z}_{\nu-1}),\,\,\,\,\,\,\forall \,\nu
$$
et donc d'après Lemme 4 et la definition de $H_\nu$ (égalité (\ref{hanu1})) pour tout $ {\bf z} =(q,a_1,a_2) \in \Lambda_\nu \setminus {\rm span} ({\bf z}_\nu)$
avec $ |{\bf z}|\ge H_\nu$ on a
\begin{equation}\label{ha}
\max_{j=1,2}|q\xi_j - a_j | \ge \varphi (H_\nu)\ge \varphi(|{\bf z}|) \ge \varphi (q).
\end{equation} 

Dans le reste de l'article, nous montrons que le vecteur $\xi$ construit satisfait la conclusion de Théorème 1.

 {\bf 4. Inégalités.}

Nous considérons l'intervalle
 \begin{equation}\label{ii}
I_\nu =
\left[\phi\left( \frac{1}{16q_\nu\rho_{\nu-1}}\right),
\frac{1}{64\rho_{\nu-1}\rho_\nu}\right].
\end{equation}

{\bf Lemme  5.} \,\,{\it
Si $ z = (q,a_1,a_2)\in \mathbb{Z}^3$ est linéairement indépendant avec
${\bf z}_{\nu-1} $ et ${\bf z}_\nu$ et
\begin{equation}\label{quu}
q\in I_\nu
\end{equation}
alors
\begin{equation}\label{relax}
\max_{j=1,2}||q\xi_j || \ge \varphi (q).
\end{equation}
}

Preuve.
 
Soit
$
\rho = \max_{j=1,2}||q\xi_j ||.
$
D’après la condition d'indépendance, nous avons
$$
0\neq
\left|
\begin{array}{ccc}
q& a_1&a_2\cr
q_{\nu-1}& a_{1,\nu-1}&a_{2,\nu-1}\cr
q_{\nu}& a_{1,\nu}&a_{2,\nu}
\end{array}
\right|
=
\left|
\begin{array}{ccc}
q& a_1-q\xi_1&a_2-q\xi_1\cr
q_{\nu-1}& a_{1,\nu-1}-q_{\nu-1}\xi_1&a_{2,\nu-1}-q_{\nu-1}\xi_2\cr
q_{\nu}& a_{1,\nu}-q_{\nu}\xi_1&a_{2,\nu}-q_{\nu2}\xi_2
\end{array}
\right|.
$$
Ainsi, selon  (\ref{uro}) nous avons
 $$
1\le  
32q\rho_{\nu-1}\rho_\nu+ 8q_\nu \rho\rho_{\nu-1} .
$$
De la borne supérieure qui découle de (\ref{quu}), nous avons
$$
\frac{1}{2}\le 8 q_\nu \rho\rho_{\nu-1}
.
$$
Ainsi
$$
\max_{j=1,2}||q\xi_j ||= \rho \ge
\frac{1}{16q_\nu\rho_{\nu-1}}\ge \varphi (q)
$$
(dans la dernière inégalité, nous utilisons la borne inférieure qui découle de  (\ref{quu})).$\Box$

{\bf 5. Preuve du Théorème 1.}

Notons que  (\ref{ce})
montre que l'union $\bigcup_\nu I_\nu$ couvre un certain rayon $ [I, +\infty)$.  Alors, d'après Lemme 5, si $q$ est assez grand et le point $(q, a_1, a_2)$ est indépendent avec deux points quelconques $ {\bf z}_{\nu-1}, {\bf z}_\nu, \nu =1,2,3,...$ alors
nous avons (\ref{relax}) et ces points ne sont pas d'intérêt pour notre propos. Donc, si nous avons une matrice unimodulaire entier 
\begin{equation}\label{maa}
\left(
\begin{array}{ccc}
q& q''& q'''\cr
a_1'&a_1''&a_1'''\cr
a_2'&a_2''&a_2'''
\end{array}
\right)
\end{equation}
avec  $ \min \{ |q'|. |q''|. |q'''|\} $ assez grand et
\begin{equation}\label{coco}
\max \left\{ \frac{\max_{j=1,2} |q'\xi_j - a_j'|}{\varphi (q')},
 \frac{\max_{j=1,2} |q''\xi_j - a_j''|}{\varphi (q'')},
 \frac{\max_{j=1,2} |q'''\xi_j - a_j'''|}{\varphi (q''')}\right\}
\le 1,
\end{equation}
alors
pour certains $\nu', \nu'', \nu'''$ on a
$$
{\bf z}' = (q', a_1',a_2') \in \Lambda_{\nu'},\,\,\,
{\bf z}'' = (q'', a_1'',a_2'') \in \Lambda_{\nu''},\,\,\,
{\bf z}''' = (q''', a_1''',a_2''') \in \Lambda_{\nu'''} .
$$
Nous prenons $\nu', \nu'', \nu'''$ être les quantités minimales à satisfaire cette propriété. Ainsi
$$
{\bf z}' = (q', a_1',a_2') \in \Lambda_{\nu'}\setminus \Lambda_{\nu'-1} ,\,\,\,
{\bf z}'' = (q'', a_1'',a_2'') \in \Lambda_{\nu''}\setminus \Lambda_{\nu''-1} ,\,\,\,
{\bf z}''' = (q''', a_1''',a_2''') \in \Lambda_{\nu'''}\setminus \Lambda_{\nu'''-1} .
$$
Nous supposons que $\nu'=\min \{ \nu', \nu'', \nu'''\} < \max \{ \nu', \nu'', \nu'''\} = \nu'''$ (sinon la matrice (\ref{maa}) a déterminant zéro, et ce n'est pas possible).

Considérons le cas  {\bf (A)} où $\nu'\le \nu''<\nu'''$.
Alors $ {\bf z}'''\in L^-({\bf z}'.{\bf z}'')\cup  L^+({\bf z}'.{\bf z}'')$ par l'unimodularité de la matrice (\ref{maa}).
Cependant $ {\bf z}''' \in \Lambda_{\nu'''}$ et $ {\bf z}''' \not\in {\rm span} ({\bf z}_{\nu'''-1})$.
Ainsi par {\bf (iv)} nous avons $ |{\bf z}''' |\ge T_{\nu'''-1}$.
Par (\ref{dee})
nous voyons que
 $
{\rm angle} (\xi, {\bf z}_{\nu'''-1})<\delta_{\nu'''-1}
.
$

Supposons que $ ({\bf z}'.{\bf z}'')\in\hbox{\got C}_{\nu'''-1}$.
Alors par  Corollaire 2 nous avons (\ref{result}) pour le point ${\bf z}'''$, c'est-à-dire
$$
\frac{\max_{j=1,2} |q'''\xi_j - a_j'''|}{\varphi (q''')}\ge 1
.
$$
Nous avons donc une contradiction avec (\ref{coco}).

Supposons que $ ({\bf z}'.{\bf z}'')\not\in\hbox{\got C}_{\nu'''-1}$.
Ensuite, par la définition de $\hbox{\got C}_{\nu'''-1}$
soit $ |{\bf z }_{\nu'} |\ge  H_{\nu'}$ ou $ |{\bf z }_{\nu''} |\ge  H_{\nu''}$
Ainsi par (\ref{ha})  avec $ \nu$ égal à  $\nu'$ ou $\nu ''$, nous avons
$$
\max \left\{ \frac{\max_{j=1,2} |q'\xi_j - a_j'|}{\varphi (q')},
 \frac{\max_{j=1,2} |q''\xi_j - a_j''|}{\varphi (q'') }\right\}
\ge 1.
$$
Nous avons une contradiction avec (\ref{coco}) à nouveau.

Donc, le cas {\bf (A)} n'est pas possible.

Maintenant, nous considérons le cas {\bf (B}) où $\nu'<\nu''=\nu'''$. Nous avons
$$
{\bf z}' \in \Lambda_{\nu'},\,\,\,\,\,\, {\bf z}'',{\bf z}''' \in \Lambda_{\nu''}\setminus\Lambda_{\nu''-1}.
$$
Comme la matrice (\ref{maa})   est unimodulaire,  le couple $({\bf z}'',{\bf z}''')$ forme une base de $\Lambda_{\nu''}$.
Ainsi
$$
{\bf z}' \in L^-({\bf z}_{\nu''}, {\bf z}_{\nu''-1})\cup L^+({\bf z}_{\nu''}, {\bf z}_{\nu''-1}).
$$
Mais alors
$$
{\bf z}_{\nu''} \in L^-({\bf z}_{\nu''-1}, {\bf z}')\cup L^+({\bf z}_{\nu''-1}, {\bf z}').
$$
Comme $ {\bf z}_{\nu ''} \in \Lambda_{\nu''}$,
par {\bf (v)} avec $ \nu = \nu''-1$ nous voyons que $ {\bf z}' \not\in \hbox{\got E}_{\nu''-1}.$ Cela signifie que
$|{\bf z}'|\ge H_{\nu'}$.
Donc par (\ref{ha}) nous avons
$$
\frac{\max_{j=1,2} |q'\xi_j - a_j'|}{\varphi (q')}
\ge 1.
$$
Nous avons une contradiction avec (\ref{coco}) à nouveau.

Donc, le cas {\bf (B)} n'est pas possible.

Donc, dans le cas où
$ \min \{ |q'|. |q''|. |q'''|\} \ge q_0 (\varphi ) $ nous avons
$$
\max \left\{ \frac{\max_{j=1,2} |q'\xi_j - a_j'|}{\varphi (q')},
 \frac{\max_{j=1,2} |q''\xi_j - a_j''|}{\varphi (q'')},
 \frac{\max_{j=1,2} |q'''\xi_j - a_j'''|}{\varphi (q''')}\right\}
\ge 1.$$
Maintenant, nous prenons $\varepsilon = \varepsilon (\varphi)>0$  assez peitit pour assurer l'inégalité nécessaire pour
les vecteurs avec $ q\le q_0 (\varphi)$.
Théorème 1 est prouvé. $\Box$

{\bf 5. Preuve du Théorème 2.}

Dans cette section, nous considérons les vecteurs meilleures approximations  ${\bf z}_\nu = (q_\nu, a_{1,\nu}, a_{2,\nu})\in \mathbb{Z}^3$
dans le sens de l'approximation simultanée de $\xi_1,\xi_2$ (voir \cite{M2}).
Nous savons que
$$
\max\left\{
\max_{j=1,2} ||q_\nu\xi_j - a_{j,\nu}||,
\max_{j=1,2} ||q_{\nu+1}\xi_j - a_{j,\nu+1}||\right\}
\le q_{\nu+1}^{-1/2} \to 0,\,\,\,\, \nu \to \infty.
$$
En outre, le couple ${\bf z}_\nu, {\bf z}_{\nu+1}$ peut être étendue à une base de $\mathbb{Z}^3$.
Le parallélogramme 
$$
\Pi_\nu =\{ {\bf x}\in \mathbb{R}^3:\,\, {\bf x}= \lambda {\bf z}_\nu+\mu{\bf z}_{\nu+1},\, \lambda , \mu \in [0,1)\}.
$$
forme un domaine fondamental du réseau de dimension deux $L({\bf z}_\nu, {\bf z}_{\nu+1})\cap \mathbb{Z}^3$.
Considérons le sous-espace affine $L^+({\bf z}_\nu, {\bf z}_{\nu+1})$ et l'orthogonal
projection $\Pi_\nu^*$ de $\Pi$ sur $L^+({\bf z}_\nu, {\bf z}_{\nu+1})$.
Alors
$$
\Pi_\nu^* = \Pi_\nu +{\bf e}_\nu,
$$
où le vecteur ${\bf e}_\nu\in \mathbb{R}^3$
a la longueur 
$$
|{\bf e}_\nu| = ({\rm det}  ( L({\bf z}_\nu, {\bf z}_{\nu+1})\cap \mathbb{Z}^3))^{-1} \asymp( q_{\nu+1} \cdot \max_{j=1,2} ||q_\nu\xi_j - a_{j,\nu}||)^{-1}
\to 0,\,\,\,\, \nu \to\infty
$$
(la dernière déclaration ici est Théorème 17 de \cite{M2} qui est une minuscule généralisation de Satz 9  de \cite{TB}).
Comme $\Pi_\nu^*$ est un domaine fondamental pour $ L({\bf z}_\nu, {\bf z}_{\nu+1})\cap \mathbb{Z}^3$, il existe un point entier
${\bf z}^* = (q^*, a_1^*, a_2^*)
\in \Pi_\nu^*$.
Nous voyons que
$$
\max_{j=1,2} ||q^*\xi_j - a_{j}^*|| = O(
\max_{j=1,2} ||q_\nu\xi_j - a_{j,\nu}||+
\max_{j=1,2} ||q_{\nu+1}\xi_j - a_{j,\nu+1}||+ |{\bf e}_\nu|) \to 0,\,\,\,\nu \to \infty.
$$
Mais la matrice
$$
\left(
\begin{array}{ccc}
q_\nu&q_{\nu+1}&q^*\cr
a_{1,\nu}& a_{1,\nu+1}&a_{1}^*\cr
a_{2,\nu}& a_{2,\nu+1}&a_{2}^*
\end{array}
\right)
$$
est unimodulaire. Théorème 2 est prouvé. $\Box$

\end{document}